\documentclass[11pt]{article}
\usepackage{indentfirst}
\usepackage{amsfonts}
\usepackage{graphicx}
\usepackage{url}
\usepackage{algorithm}
\usepackage{algpseudocode}
\usepackage{listings}
\usepackage{caption}
\usepackage{subcaption}
\newtheorem{teor}{Theorem}[section]
\newtheorem{exemplu}{Example}[section]

\author{E. Scheiber\thanks{e-mail: scheiber@unitbv.ro}} 
\title{Adjoint System in the Shooting Method to Solve Boundary Value Problems} 
\date{}
\begin{document}
\maketitle

\begin{abstract}
The shooting method is used to solve a boundary value problem with
separated and explicit constraints. To obtain
approximations of an unknown initial values there are considered
arguments based on the adjoint differential system attached to the given differential 
system. Finally the Newton-Kantorovich iterations are regained.


\end{abstract}

\section{Introduction}
The shooting method to solve a boundary value problem (BVP) is presented in many textbooks \cite{3}, \cite{6}, \cite{7}, \cite{2}.
In this note we consider the case of separated and explicit boundary conditions. For some variables
the initial value is given and for the rest its are unknown. 
In essence the scope of the shooting method is to compute the unknown initial values in order to satisfy the ending conditions. 

To use / implement the method two ingredients are necessary:
\begin{itemize}
\item
A method / routine to integrate an initial value problem;
\item
A method / routine to solve an algebraic system. For a nonlinear algebraic system the Newton-Kantorovich method is
widely used.
\end{itemize}
Practical aspects of implementing the shooting method are given in \cite{4,5}. 

In this note we point out a connection with an adjoint system  
attached to the differential system of the BVP. Trying to compute the unknown initial values
using adjoint functions we regain the Newton-Kantorovich iterations.

The paper is organized as follows. After defining the BVP in Section 2, the shooting method is
revisited in Section 3. In Section 4 the connection with adjoined system is detailed. Finally, in Section 5 
numerical results for several examples are given.

\section{The boundary value problem with explicit boundary constraints}
The considered BVP is defined by a system of ordinary differential equations 
\begin{equation}\label{bvp1}
\dot{x}(t) = f(t,x(t)),\qquad  t\in[0,T],
\end{equation}
and the boundary conditions
\begin{equation}\label{bvp2}
\begin{array}{lclcl}
x_i(0)=y_{0,i},  & \qquad & i\in\mathcal{I}_0\subset\{1,2,\ldots,n\},& \quad & |\mathcal{I}_0|=m\\
x_j(T)=y_{T,j}, &   & j\in\mathcal{I}_T\subset\{1,2,\ldots,n\},& \quad & |\mathcal{I}_T|=n-m.
\end{array}
\end{equation}
where 
$f:[0,T]\times\mathbb{R}^n\rightarrow\mathbb{R}^n, \mathcal{I}_0, \mathcal{I}_T,
y_0\in\mathbb{R}^m, y_T\in\mathbb{R}^{n-m}$ are given. The function $f(t,x)$ is twice continuously differentiable.
By $|\mathcal{I}_0|,|\mathcal{I}_T|$ we noted the number of elements of the corresponding sets.
Let be $\mathcal{I}_0^c=\{1,2,\ldots,n\}\setminus\mathcal{I}_0$,  the set of index for which the initial values are unknown.

We will use the matrices $P_T\in M_{n-m,n}(\{0,1\}), P_0^c\in M_{n-m,n}(\{0,1\})$ with the following definitions:

\noindent
If $\mathcal{I}_T=\{k_1,k_2,\ldots,k_{n-m}\}$ with $k_1<k_2<\ldots<k_{n-m}$ then
$$
(P_T)_{i,j}=\left\{\begin{array}{lclcl}
1 &\mathrm{if} & j=k_i, &\qquad & i\in\{1,2,\ldots,n-m\}\\
0 & & \mathrm{otherwise} & & j\in\{1,2,\ldots,n\}
\end{array}\right.
$$
and $P_0^c$ is defined analogously, corresponding to the set $\mathcal{I}_o^c.$

\noindent
We suppose that the BVP has at least a solution.

\vspace*{0.3cm}%\noindent
For simplicity let be 
%\begin{equation}\label{Hyp}
%\begin{array}{l}
$\mathcal{I}_0=\{1,2,\ldots,m\}$ and
$\mathcal{I}_T=\{j_1,j_2,\ldots,j_{n-m}\}.$
%\end{array}
%\end{equation}
Then the boundary conditions are
$$
\begin{array}{lclcl}
x_i(0) &=&y_{0,i},& \qquad & i\in\{1,2,\ldots,m\}\\
x_{j_k}(T)&=&y_{T,k}, & & k\in\{1,2,\ldots,n-m\}
\end{array}.
$$
Because $\mathcal{I}_0^c=\{m+1,\ldots,n\}$ it follows that $P_0^c=[0_{n-m,m}\ I_{n-m}].$

\section{Recall of the shooting method}

We shall follow J. Singh's presentations \cite{4,5}.

If the initial value of (\ref{bvp1}) is 
\begin{equation}\label{bvpet24}
x_0(c)=\left[\begin{array}{l}y_0\\c\end{array}\right]\qquad\Leftrightarrow\qquad P_0^cx_0(c)=c,
\end{equation}
with $c\in\mathbb{R}^{n-m},$ 
then the boundary conditions in $T$  may be written as
%$c$ must be determined to satisfy the boundary conditions in $T$ 
\begin{equation}\label{bvpet12}
F(c)=P_T x(T;x_0(c))-y_T=\left(\begin{array}{l}
x_{j_1}(T;x_0(c))-y_{T,1}\\
\vdots\\
x_{j_{n-m}}(T;x_0(c))-y_{T,n-m}
\end{array}\right)=0.
\end{equation}
(\ref{bvpet12}) represents an algebraic system of $n-m$ equations where
the components of $c=(c_1,\ldots,c_{n-m})^{\mathrm{T}}$ are the unknowns. 
The symbol $^{\mathrm{T}}$ marks the transpose of a matrix.

The Newton-Kantorovich methods may be used to solve this system. Then the following recurrence formula
is applied:
\begin{equation}\label{bvpet18}
c^{k+1}=c^k-(F'(c^k))^{-1}F(c^k),\qquad k=0,1,2,\ldots.
\end{equation}

In order to compute $F(c^k),$ the initial value problem (\ref{bvp1})-(\ref{bvpet24}), with $c=c^k$, must be integrated.
The computation of the matrix 
$$
F'(c)=P_T\frac{\partial x(T;x_0(c))}{\partial c}=
\left(\begin{array}{ccc}
\frac{\partial x_{j_1}(t;x_0(c))}{\partial c_1} & \ldots & \frac{\partial x_{j_1}(t;x_0(c))}{\partial c_{n-m}}\\
\vdots &  & \vdots \\
\frac{\partial x_{j_{n-m}}(t;x_0(c))}{\partial c_1} & \ldots & \frac{\partial x_{j_{n-m}}(t;x_0(c))}{\partial c_{n-m}}
\end{array}\right)
$$ 
is detailed below.

From
$$
\frac{\mathrm{d}}{\mathrm{d}t}\frac{\partial}{\partial c_j}x_i(t;x_0(c))=
\frac{\partial}{\partial c_j}\frac{\mathrm{d}}{\mathrm{d}t}x_i(t,x_0(c))=
$$
$$
=\frac{\partial}{\partial c_j} f_i(t,x(t;x_0(c))=
\sum_{k=1}^n \frac{\partial f_i(t,x(t;x_0(c))}{\partial x_k} \frac{\partial x_k(t;x_0(c))}{\partial c_j}
$$
it results that
\begin{equation}\label{bvpet14}
\frac{\mathrm{d}}{\mathrm{d}t}\underbrace{\left(\begin{array}{ccc}
\frac{\partial x_{1}(t;x_0(c))}{\partial c_1} & \ldots & \frac{\partial x_{1}(t;x_0(c))}{\partial c_{n-m}}\\
\vdots &  & \vdots \\
\frac{\partial x_{n}(t;x_0(c))}{\partial c_1} & \ldots & \frac{\partial x_{n}(t;x_0(c))}{\partial c_{n-m}}\\
\end{array}\right)}_{\frac{\partial x(t,x_0(c))}{\partial c}}=
\end{equation} 
$$
=\underbrace{\left(\begin{array}{ccc}
\frac{\partial f_{1}(t,x)}{\partial x_1} & \ldots & \frac{\partial f_{1}(t,x)}{\partial x_n}\\
\vdots &  & \vdots \\
\frac{\partial f_{n}(t,x)}{\partial x_1} & \ldots & \frac{\partial f_{n}(t,x)}{\partial x_n}\\
\end{array}\right)}_{f'(t,x(t;x_0(c)))}
\underbrace{\left(\begin{array}{ccc}
\frac{\partial x_1(t;x_0(c))}{\partial c_1} & \ldots & \frac{\partial x_1(t;x_0(c))}{\partial c_{n-m}}\\
\vdots &  & \vdots \\
\frac{\partial x_n(t;x_0(c))}{\partial c_1} & \ldots & \frac{\partial x_n(t;x_0(c))}{\partial c_{n-m}}\\
\end{array}\right)}_{\frac{\partial x(t,x_0(c))}{\partial c}}
$$
The initial conditions for this differential system are
\begin{equation}\label{bvpet17}
\frac{\partial x(0;x_0(c))}{\partial c_i}=e_{m+i}\qquad i\in\{1,2,\ldots,n-m\},
\end{equation}
where $e_{m+i}$ is the notation of the vector of the canonical base from $\mathrm{R}^n.$\\
The two initial value problems are integrated simultaneously %, in the sense that
\begin{equation}\label{bvpet23}
\left[\begin{array}{c}
\dot{x}(t)\\
\\
\frac{\mathrm{d}}{\mathrm{d}t}\frac{\partial x(t;x_0(c))}{\partial c_i}
\end{array}\right]=\left[\begin{array}{c}
f(t,x(t))\\
\\
f'(t,x(t;x_0(c))\frac{\partial x(t,x_0(c))}{\partial c_i}
\end{array}\right],
\end{equation}
with the initial conditions
$$ 
\left[\begin{array}{c}
x(0,x_0(c)) \\
\\
\frac{\partial x(0;x_0(c))}{\partial c_i}
\end{array}\right]=
\left[\begin{array}{c}\left[\begin{array}{c}y_0\\c\end{array}\right]\\e_{m+i}\end{array}\right],\qquad i\in\{1,2,\ldots,n-m\}.
$$
The last $n$ rows of a column of the solution in $T$ is a column of $\frac{\partial x(t,x_0(c))}{\partial c_i}.$

\noindent
In total $n-m$ initial value problems must be solved, for each $c_i.$ 

From the equality
$$
\frac{\partial x(t;x_0(c))}{\partial c}=\frac{\partial x(t;x_0(c))}{\partial x_0} (P_0^c)^{\mathrm{T}}
$$
it results
\begin{equation}\label{bvpet20}
F'(c)=P_T \frac{\partial x(t;x_0(c))}{\partial x_0}( P_0^c)^{\mathrm{T}}.
\end{equation}

\section{The connection with adjoint system}

Let  the ordinary differential system (\ref{bvp1}) be given. The following system
\begin{equation}\label{bvp4}
\dot{p}(t)=-f_x^{'\mathrm{T}}(t,x(t))p(t)
\end{equation} 
is called the adjoint differential system to (\ref{bvp1}) \cite{1}.
Here we have used the notation $f_x'(t,x)=\frac{\partial f(t,x)}{\partial x}.$
A solution of this system is an adjoint function. 

We denote by $p^j(t)=(p_1^j,\ldots,p_n^j)^{\mathrm{T}}$ the adjoint function satisfying the initial condition $p^j(T)=e_j$ and
$$
\mathrm{P}(t)=[p^1(t)\   p^2(t)\ \ldots\  p^n(t)].
$$

There is a connection between $\mathrm{P}(0)$ and $\frac{\partial x(T;x_0))}{\partial x_0}:$
\begin{teor}\label{bvpt1}
The following equality holds
$$
\frac{\partial x(T;x_0)}{\partial x_0}=\mathrm{P}(0)^{\mathrm{T}}.
$$
\end{teor}

\vspace*{0.3cm}\noindent\textbf{Proof.}
%The following equalities hold
We have the equalities
$$
<p^j(T), \frac{\partial x(T;x_0)}{\partial x_{0,k}} >-<p^j(0), \frac{\partial x(0;x_0)}{\partial x_{0,k}} >=
$$
$$
=\int_0^T \frac{\mathrm{d}}{\mathrm{d}t}<p^j(t),\frac{\partial x(t;x_0)}{\partial x_{0,k}}>\mathrm{d}t=
$$
$$
=\int_0^T\left(<\dot{p}^j(t),\frac{\partial x(t;x_0)}{\partial x_{0,k}}>+
<p^j(t),\frac{\mathrm{d}}{\mathrm{d}t} \frac{\partial x(t;x_0)}{\partial x_{0,k}}>\right)\mathrm{d}t=
$$
$$
=\int_0^T\left(<-f_x^{'\mathrm{T}}(t,x(t;x_0))p^j(t),\frac{\partial x(t;x_0)}{\partial x_{0,k}}>+\right.
$$
$$
\left.+<p^j(t),f_x'(t,x(t;x_0)) \frac{\partial x(t;x_0)}{\partial x_{0,k}}>\right)\mathrm{d}t=0.
$$
Taking into account the initial conditions, $p^j(T)=e_j,  \frac{\partial x(0;x_0)}{\partial x_{0,k}}=e_k,$ it results
\begin{equation}\label{bvpet21}
\frac{\partial x_j(T;x_0)}{\partial x_{0,k}}=p^j_k(0).
\end{equation}
Then it follows that
$$
\frac{\partial x(T;x_0))}{\partial x_0}=
\left(\begin{array}{ccc}
\frac{\partial x_{1}(t;x_0)}{\partial x_{0,1}} & \ldots & \frac{\partial x_{1}(t;x_0)}{\partial x_{0,n}}\\
\vdots &  & \vdots \\
\frac{\partial x_{n}(t;x_0)}{\partial x_{0,1}} & \ldots & \frac{\partial x_{n}(t;x_0)}{\partial x_{0,n}}
\end{array}\right)=
\left(\begin{array}{ccc}
p^1_1(0) & \ldots & p^1_n(0) \\
\vdots & & \vdots \\
p^n_1(0) & \ldots & p^n_n(0) 
\end{array}\right)=\mathrm{P}(0)^{\mathrm{T}}.
$$
\rule{5pt}{5pt}

\vspace*{0.3cm}
Using adjoint functions we compute the unknown initial values of the BVP (\ref{bvp1})-(\ref{bvp2}).
Let be $x^*(t), x(t)$ two functions verifying (\ref{bvp1}) and $\Delta_x=x^*-x.$ Then
$$
\dot{\Delta}_x(t)=\dot{x}^*(t)-\dot{x}(t)=
f(t,x^*(t))-f(t,x(t))=f_x'(t,x(t))\Delta_x(t)+O(\|\Delta_x(t)\|^2).
$$
Neglecting the last term, let $\delta_x(t)$ be such that
\begin{equation}\label{bvp5}
\dot{\delta}_x(t)=f_x'(t,x(t))\delta_x(t).
\end{equation}

If $p(t)$ is an adjoint function, with similar computations as in the proof of the above theorem we  obtain
\begin{equation}\label{bvp6}
<p(T),\delta_x(T)>-<p(0),\delta_x(0))>=0.
\end{equation}

If $x^*$ is the solution of the BVP (\ref{bvp1})-(\ref{bvp2}) and $x$
is an approximation of $x^*$ then we may suppose that 
$x^*_i(0)=x_i(0)=y_{0,i}$ and thus $\delta_{x,i}(0)=\Delta_{x,i}(0)=0,$
for any $i\in\mathcal{I}_0=\{1,2,\ldots,m\}.$
The components of $\delta_x$ are denoted by $\delta_{x,i}, i\in\{1,2,\ldots,n\}.$

For any $i\in\mathcal{I}_T$
%\footnote{Due to the hypothesis (\ref{Hyp}) $\mathcal{J}=\{m+1,\ldots,n\}.$}
we set the adjoint function $p:=p^{i}$ 
%the adjoint function with the initial condition $p^i(T)=e_i.$ %$ (given in $T$),
%=(p_1^i,\ldots,p_n^i)^{\mathrm{T}}
%where $e_i$ is the notation of the vector of the canonical base from $\mathrm{R}^n.$ 
and from (\ref{bvp6}) we deduce
$$
\delta_{x,i}(T)-\sum_{j=1}^{n-m}p^i_{m+j}(0)\delta_{x,m+j}(0)=0,\qquad \forall\  i\in\mathcal{I}_T,
$$
or
\begin{equation}\label{bvp7}
\sum_{j=1}^{n-m}p^i_{m+j}(0)\delta_{x,m+j}(0)=y_{T,i}-x_i(T),\qquad \forall\  i\in\mathcal{I}_T,
\end{equation}
and in matrix form 
$$
\left(\begin{array}{lcl}
p^{j_1}_{m+1} & \ldots & p^{j_1}_n\\
\vdots && \vdots \\
p^{j_{n-m}}_{m+1} & \ldots & p^{j_{n-m}}_n
\end{array}\right)
\left(\begin{array}{c}
\delta_{x,m+1}(0)\\
\vdots\\
\delta_{x,n}(0)
\end{array}\right)=
\left(\begin{array}{c}
y_{T,j_1}-x_{j_1}(T)\\
\vdots\\
y_{T,j_{n-m}}-x_{j_{n-m}}(T)
\end{array}\right).
$$

We have obtained an algebraic linear system with the unknowns $\delta_{x,m+j}(0),\\ j\in\{1,\ldots,n-m\}.$

We interpret (\ref{bvp7}) as the needed correction of the initial conditions of $x$ in order to satisfy the boundary conditions in $T.$

Due to the neglection made to obtain (\ref{bvp6}) we suppose that $x_{m+j}(0)+\delta_{x,m+j}(0),
j\in\{1,\ldots,n-m\},$ are better approximations for the unknown initial conditions and an iterative scheme
must be taken into account. If $x(t)=x(t;x_0(c^k))$ and $\delta_{x,m+j}(0)=c^{k+1}_j-c^k_j,\ j\in\{1,2,\ldots,n-m\}$
then (\ref{bvp7}) is rewritten as
$$
P_T\mathrm{P}(T)^{\mathrm{T}}(P_0^c)^{\mathrm{T}}(c^{k+1}-c^k)=-F(c^k).
$$
Taking into account the Theorem \ref{bvpt1} and the equality (\ref{bvpet20}) we regain the
Newton-Kantorovich iterations.

\section{Examples}
The choice of the unknown initial values is very important mainly for the success of the numerical integration.

\begin{exemplu}\cite{5}
\end{exemplu}
\begin{eqnarray*}
y''(t) &=& 2y(t) y'(t),\qquad t\in[0,1]\\
y(0) &=& 0\\
y(1) &=& 2
\end{eqnarray*}
The solution of this BVP is $y(t)=a \cdot \mathrm{tan}(a t),$ where $a\approx 1.0768740.$ 
%The initial values of the corresponding first order differential problem are $x_0=\left(\begin{array}{c}0\\1.1596576\end{array}\right).$

\vspace*{0.3cm}
\begin{tabular}{|c|c|}
\hline\hline
Initial values & Final values \\
\hline\hline
$y(0)=0$  & $y(1)=2.0000000$ \\
$\frac{\mathrm{d}y}{\mathrm{d}t}(0)=1.1596576$ & $\frac{\mathrm{d}y}{\mathrm{d}t}(1)=5.1596576$ \\
\hline
\end{tabular}

\begin{figure}
     \centering
     \begin{subfigure}[b]{0.45\textwidth}
         \centering
         \includegraphics[width=\textwidth]{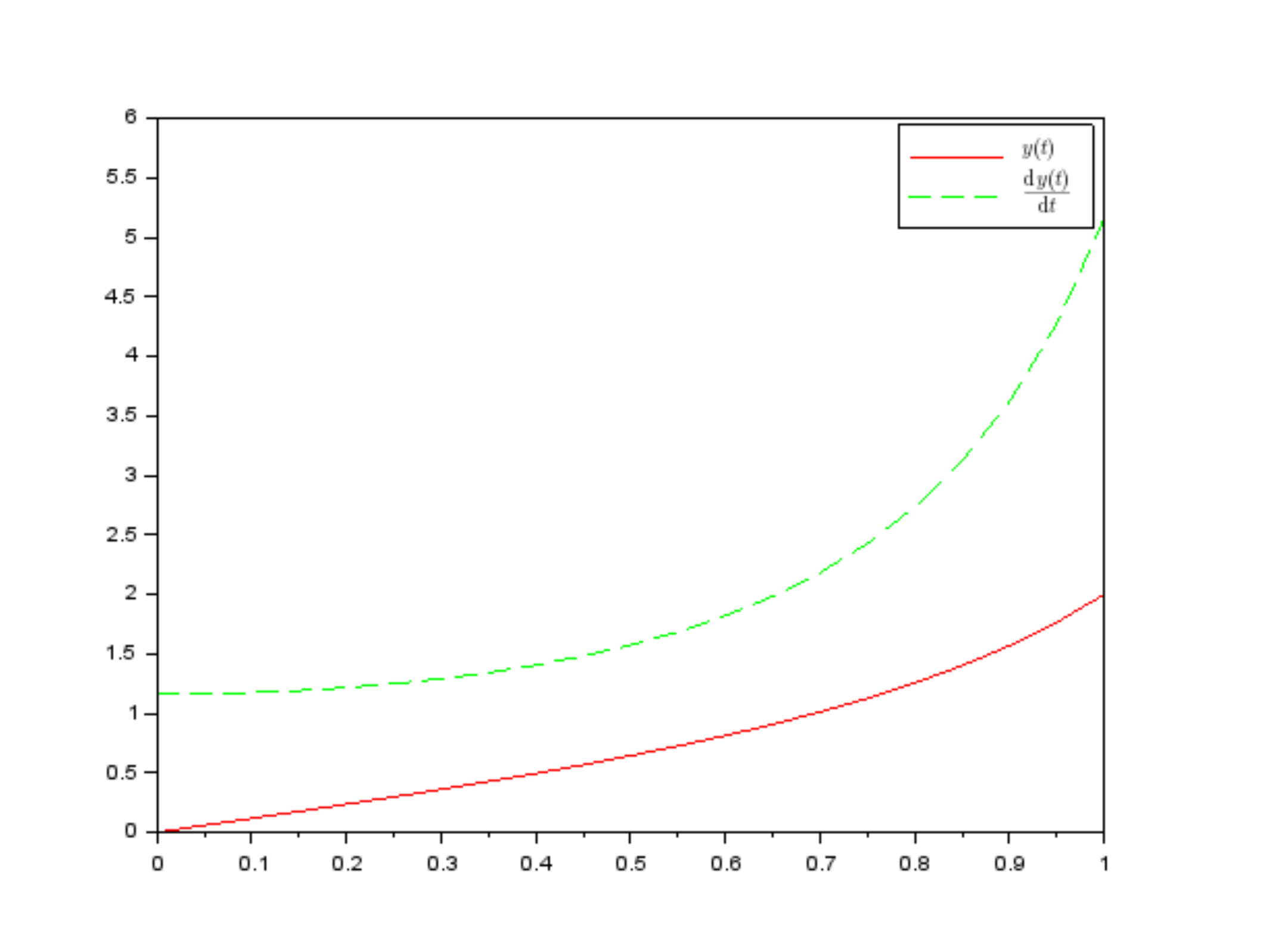}
         \caption{Example 1: Plots of $y(t)$ and $\frac{\mathrm{d}y(t)}{\mathrm{d}t}.$}
         \label{eq1}
     \end{subfigure}
     \hfill
     \begin{subfigure}[b]{0.45\textwidth}
         \centering
         \includegraphics[width=\textwidth]{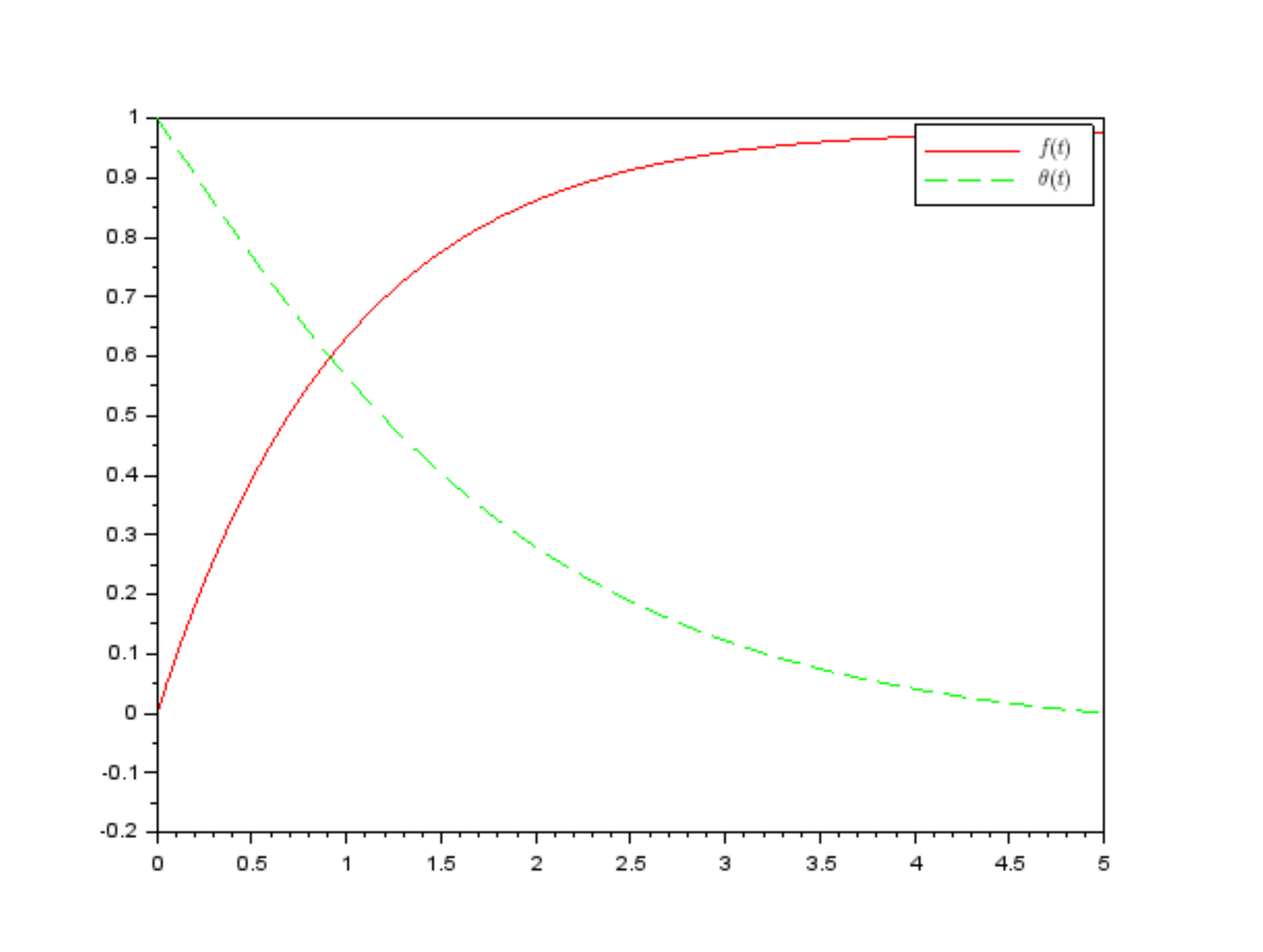}
         \caption{Example 2: Plots of $f$ and $\theta$ - case 1.}
         \label{eq2_1}
     \end{subfigure}
     \hfill
     \begin{subfigure}[b]{0.45\textwidth}
         \centering
         \includegraphics[width=\textwidth]{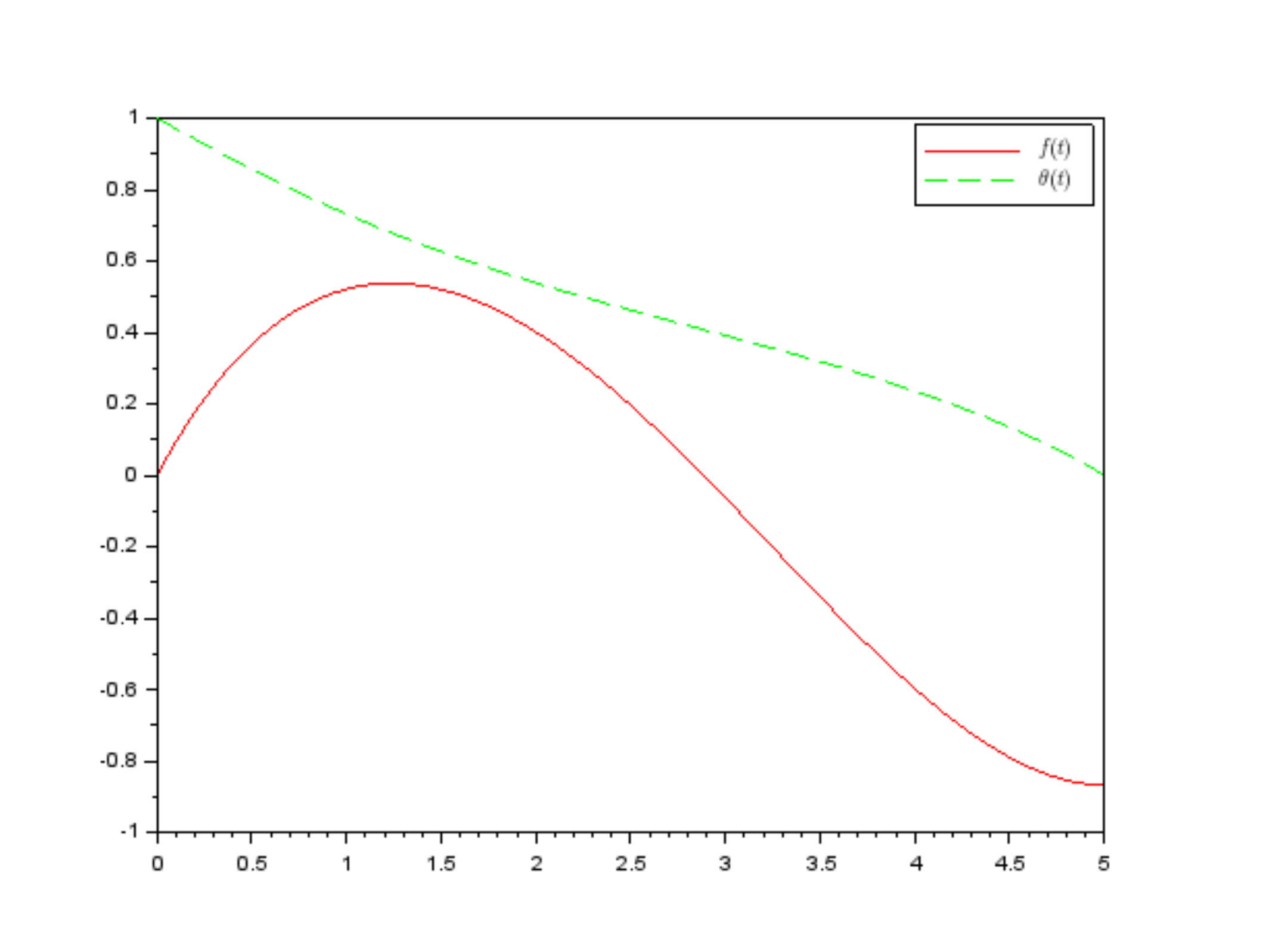}
         \caption{Example 2: Plots of $f$ and $\theta$ - case 2.}
         \label{eq2_2}
     \end{subfigure}
    \hfill
     \begin{subfigure}[b]{0.45\textwidth}
         \centering
         \includegraphics[width=\textwidth]{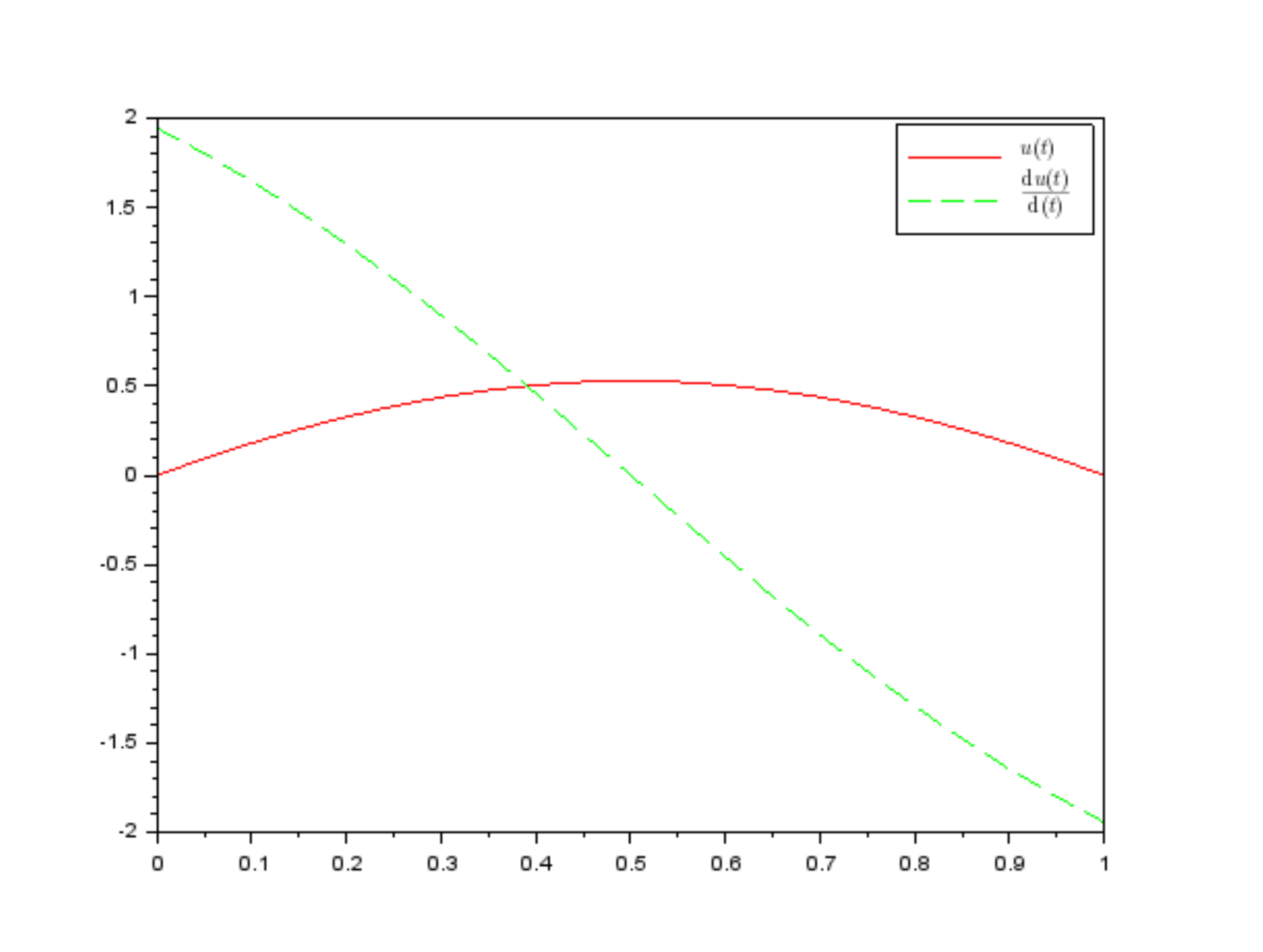}
         \caption{Example 3: Plots of $u$ and $\frac{\mathrm{d}u(t)}{\mathrm{d}t}$ - case 1.}
         \label{eq3_1}
     \end{subfigure}
     \hfill
     \begin{subfigure}[b]{0.45\textwidth}
         \centering
         \includegraphics[width=\textwidth]{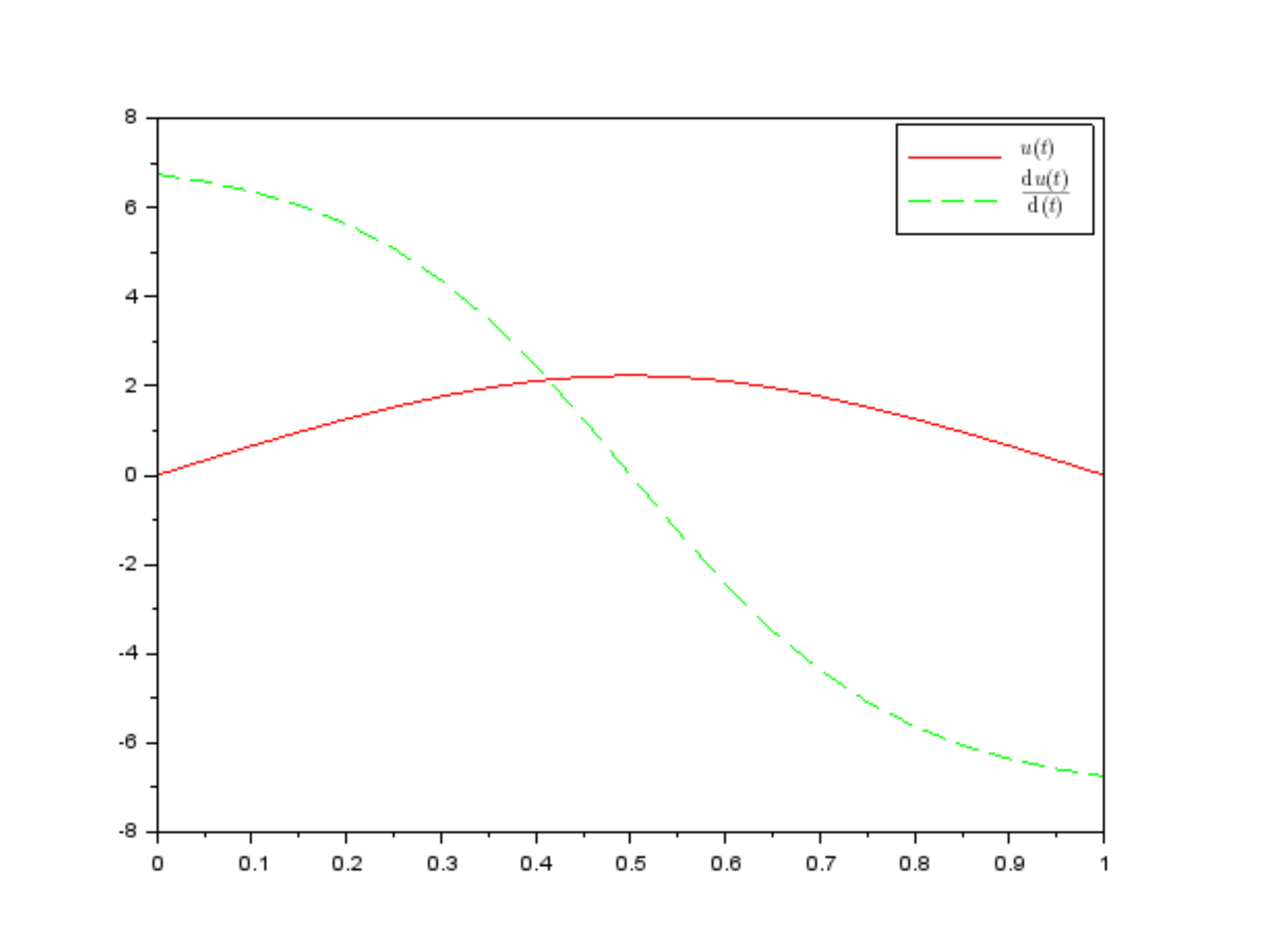}
         \caption{Example 3: Plots of $u$ and $\frac{\mathrm{d}u(t)}{\mathrm{d}t}$ - case 2.}
         \label{eq3_2}
     \end{subfigure}
     \hfill
     \begin{subfigure}[b]{0.45\textwidth}
         \centering
         \includegraphics[width=\textwidth]{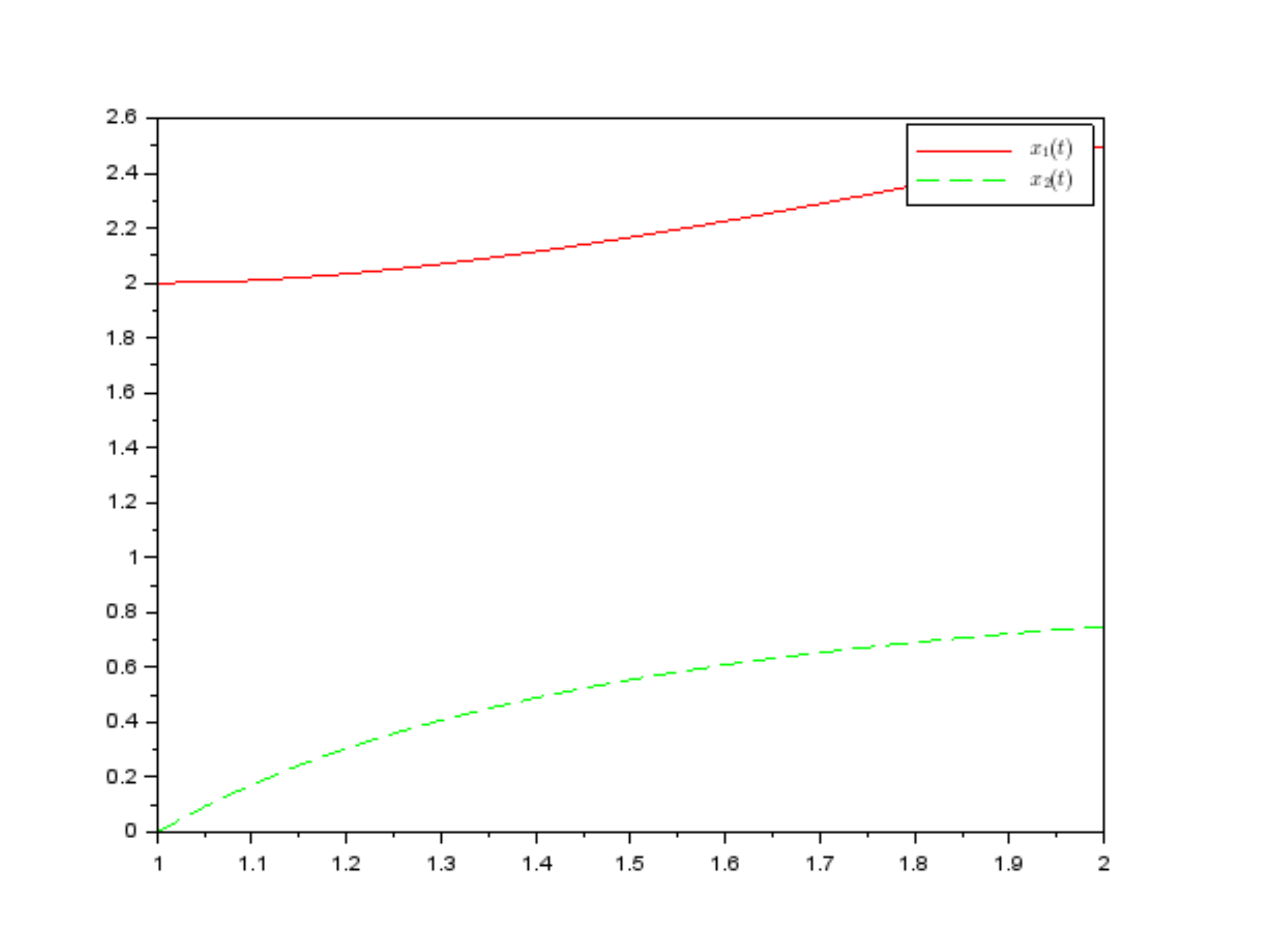}
         \caption{Example 4: Plots of $x_1(t)$ and $x_2(t).$}
         \label{eq4}
     \end{subfigure}
        \caption{Graphic representations}
%        \label{fig:three graphs}
\end{figure}

\vspace*{0.3cm}\noindent
The plot of $y(t)$ and $\frac{\mathrm{d}y(t)}{\mathrm{d}t}$ are given in Fig \ref{eq1}.

%\begin{center}
%\begin{figure}[h]
%\hspace*{2cm}
%\includegraphics[width=12cm,height=9cm,keepaspectratio]{images/Eq1.pdf}

%\caption{Example 1: Plots of $y(t)$ and $\frac{\mathrm{d}y(t)}{\mathrm{d}t}.$}\label{eq1}
%\end{figure}
%\end{center}

\begin{exemplu}\cite{5}
\end{exemplu}
$$
\begin{array}{cclcl}
%\begin{eqnarray*}
f^{(3})+f f''-(f')^2 &=& 0,&\qquad & t\in[0,5]\\
\theta''+ k \theta' f &=& 0 &&\\
f(0) &=& 0 &&\\
f'(0) &=& 1 &&\\
f'(5) &=& 0&& \\
\theta(0) &=& 1 && \\
\theta(5) &=& 0 & 
%\end{eqnarray*}
\end{array}
$$
where $k$ is the Prandtl number. We used $k=0.71.$

For different values of the starting initial values it may find different results:

\vspace*{0.3cm}
\begin{tabular}{|c|c|}
\hline\hline
Case & $c^0$    \\

\hline\hline
1 & $c=\left(\begin{array}{c}0\\0\end{array}\right),\left(\begin{array}{c}-1\\-1\end{array}\right)$ \\
\hline\hline
\end{tabular}

\vspace*{0.3cm}
\begin{tabular}{|l|l|}
\hline\hline
Initial values & Final values \\
\hline\hline
$f(0)=0$  & $f(1)=0.9740442$ \\
$f'(0)=1$ & $f'(1)=-1.044D-15$ \\
$f''(0)=-1.0013962$ & $f''(1)=-0.0072487$ \\
$\theta(0)=1$ & $\theta(1)=-7.672D-16$ \\
$\theta'(0)= -0.4755621$ & $\theta'(1)=-0.0283081$ \\
\hline
\end{tabular}

\begin{tabular}{|c|c|}
\hline\hline
Case & $c^0$    \\
\hline\hline
2 & $c=\left(\begin{array}{c}-2\\0\end{array}\right)$ \\
\hline\hline
\end{tabular}

\vspace*{0.3cm}
\begin{tabular}{|l|l|}
\hline\hline
Initial values & Final values \\
\hline\hline
$f(0)=0$  & $f(1)=-0.8678587$ \\
$f'(0)=1$ & $f'(1)=4.820D-16$ \\
$f''(0)=-1.2108404$ & $f''(1)=0.7142624$ \\
$\theta(0)=1$ & $\theta(1)=-1.419D-15$ \\
$\theta'(0)= -0.2921733$ & $\theta'(1)=-0.3115125$ \\
\hline
\end{tabular}

\vspace{0.3cm}
The plot of the trajectories in the two cases are given in Fig. \ref{eq2_1} and Fig. \ref{eq2_2}.

%\begin{center}
%\begin{figure}[h]
%\hspace*{2cm}
%\includegraphics[width=12cm,height=9cm,keepaspectratio]{images/Eq2_1.pdf}

%\caption{Example 2: Plots of $f$ and $\theta$ - case 1.}\label{eq2_1}
%\end{figure}
%\end{center}

%\begin{center}
%\begin{figure}[h]
%\hspace*{2cm}
%\includegraphics[width=12cm,height=9cm,keepaspectratio]{images/Eq2_2.pdf}

%\caption{Example 2: Plots of $f$ and $\theta$ - case 2.}\label{eq2_2}
%\end{figure}
%\end{center}

\begin{exemplu}\cite{3}\end{exemplu}
$$
\begin{array}{l}
u''(t)+e^{u(t)+1}=0  \qquad  t\in[0,1];\\
u(0)=u(1)=0 
\end{array}
$$
This problem with the Bratu type equation \cite{6} has two solutions 
$$
u(t)=-2 \ln\frac{\cosh((t-\frac{1}{2})\frac{\theta}{2})}{\cosh\frac{\theta}{4}}
$$
where $\theta$ is a solution of the equation $\theta=\sqrt{2e}\cosh\frac{\theta}{4}.$
This equation has two solutions $\theta_1\approx 3.0362318$ and $\theta_2\approx 7.1350055.$

\vspace*{0.3cm}
\begin{tabular}{|c|c|}
\hline\hline
Case & $c^0$    \\
\hline\hline
1 & $c=0$ \\
\hline\hline
\end{tabular}

\vspace*{0.3cm}
\begin{tabular}{|l|l|}
\hline\hline
Initial values & Final values \\
\hline\hline
$u(0)=0$  & $u(1)=2.072D-10$ \\
$u'(0)=1.9447725$ & $u'(1)=-1.9447725$ \\
\hline
\end{tabular}

\begin{tabular}{|c|c|}
\hline\hline
Case & $c^0$    \\
\hline\hline
2 & $c=5$ \\
\hline\hline
\end{tabular}

\vspace*{0.3cm}
\begin{tabular}{|l|l|}
\hline\hline
Initial values & Final values \\
\hline\hline
$u(0)=0$  & $u(1)=-2.196D-08$ \\
$u'(0)=6.7432737$ & $u'(1)=-6.7432738$ \\
\hline
\end{tabular}

\vspace{0.3cm}
The plot of the trajectories in the two cases are given in Fig. \ref{eq3_1} and Fig. \ref{eq3_2}.

%\begin{center}
%\begin{figure}[h]
%\hspace*{2cm}
%\includegraphics[width=12cm,height=9cm,keepaspectratio]{images/Eq3_1.pdf}

%\caption{Example 3: Plots of $u$ and $\frac{\mathrm{d}u(t)}{\mathrm{d}t}$ - case 1.}\label{eq3_1}
%\end{figure}
%\end{center}

%\begin{center}
%\begin{figure}[h]
%\hspace*{2cm}
%\includegraphics[width=12cm,height=9cm,keepaspectratio]{images/Eq3_2.pdf}

%\caption{Example 3: Plots of $u$ and $\frac{\mathrm{d}u(t)}{\mathrm{d}t}$- case 2.}\label{eq3_2}
%\end{figure}
%\end{center}

\begin{exemplu}\cite{5}\end{exemplu}

$$
\begin{array}{lcll}
\dot{x}_1(t) &=& x_2(t) & t\in[1,2], \\
\dot{x}_2(t) &=& 2 x_1^3(t)-6 x_1(t)-2 t^3 &\\
x_1(1) &=& 2 &\\
x_1(2) &=& 2.5 &
\end{array}
$$

\vspace*{0.3cm}
\begin{tabular}{|l|l|}
\hline\hline
Initial values & Final values \\
\hline\hline
$x_1(1)=2$  & $x_1(2)=2.5000000$ \\
$x_2(1)=6.363D-10$ & $x_2(2)=0.7500000$ \\
\hline
\end{tabular}

\vspace*{0.3cm}\noindent
The plot of $x_1(t)$ and $x_2(t)$ are given in Fig \ref{eq4}.

%\begin{center}
%\begin{figure}[h]
%\hspace*{2cm}
%\includegraphics[width=12cm,height=9cm,keepaspectratio]{images/Eq4.pdf}

%\caption{Example 4: Plots of $x_1(t)$ and $x_2(t).$}\label{eq4}
%\end{figure}
%\end{center}

To make the results reproducible we provide some code 
at \url{https://github.com/e-scheiber/bvp.git}. 

\end{document}